\documentclass[a4paper,reqno]{amsart}

\usepackage{amscd}
\usepackage{amssymb,latexsym}
\usepackage[centertags]{amsmath}
\usepackage{amsfonts}
\usepackage{amssymb}
\usepackage{amsthm}
\usepackage{newlfont}
\usepackage[all]{xy}
\usepackage{amscd}
\usepackage{amssymb,latexsym}
\usepackage[all]{xy}       

\def\subsection#1{\removelastskip\par\medskip\setcounter{equation}{0}
\refstepcounter{subsection}\noindent{\bf\thesubsection. #1}}

\def\g{ $\displaystyle {\frak g}$ }
\def\g01{ $\displaystyle {\frak g}={\frak g}_{\bar 0}\oplus{\frak g}_{\bar 1}$ }
\def\g0{ $\displaystyle {\frak g}_{\bar 0}$ }
\def\g1{ $\displaystyle {\frak g}_{\bar 1}$ }

 \newtheorem{teo}{Theorem}[section]

 \newtheorem{exa}[teo]{Example}
 \newtheorem{pro}[teo]{Proposition}
 
 \theoremstyle{definition}
 \newtheorem{defi}[teo]{Definition}
 \theoremstyle{remark}
 
\def\3{{\sf 3}}

\def\g{{\gamma}}

\begin{document}

\title[Quadratic Malcev superalgebras with
reductive even part]{Quadratic Malcev superalgebras with reductive
even part} \author[H. Albuquerque ]{Helena~Albuquerque}

\address{Helena~Albuquerque, Departamento
de Matem\'{a}tica, Universidade de Coimbra, Apartado 3008, 3001-454
Coimbra, PORTUGAL {\em E-mail address}: {\tt lena@mat.uc.pt}}

\thanks{The first and the second
authors  acknowledge partial financial assistance by the CMUC,
Department of Mathematics, University of Coimbra.}

\author[E. Barreiro]{Elisabete~Barreiro}

\address{Elisabete~Barreiro, Departamento de Matem{\'a}tica,
Universidade de Coimbra, Apartado 3008, 3001-454 Coimbra, PORTUGAL
{\em E-mail address}: {\tt mefb@mat.uc.pt}}{}


\author[S. Benayadi]{Sa\"{i}d~Benayadi}

\address{Sa\"{i}d~Benayadi, Laboratoire de Mathématiques et
Applications de Metz, CNRS-UMR 7122, Université Paul Verlaine, Ile
du Saulcy, 57 045 Metz cedex 1, FRANCE {\em E-mail address}: {\tt
benayadi@poncelet.univ-metz.fr}} 

\begin{abstract}
It is our goal to give an inductive description of quadratic Malcev
superalgebras with reductive even part. We  use  the notion of
double extension of Malcev superalgebras presented by H. Albuquerque
and S. Benayadi in \cite{AB2004} and transfer to Malcev
superalgebras the concept of generalized double extension given in
\cite{BBB} for Lie superalgebras.
\end{abstract}

\maketitle

\bigskip

\section*{Introduction}

In this work all the quadratic Malcev superalgebras considered are
finite dimensional, defined over an algebraically closed field of
characteristic 0. The description of Malcev superalgebras with
reductive even part and the action of the even part on the odd part
is completely reducible was reduced in \cite{Alb1993} to the
classification of Lie superalgebras satisfying the same conditions.
This problem was solved by Alberto Elduque in \cite{Eld1996}, inside
the description of the structure of semisimple Lie superalgebras
using the notions of nice extension and elementary extension.

A quadratic Malcev superalgebra is a Malcev superalgebra where is
defined an  even non-degenerate invariant supersymmetric bilinear
form. Quadratic Lie algebras were classified by A. Medina and P.
Revoy in \cite{MR1985} using the notion of double extension of Lie
algebras. This notion was extended in \cite {BB1999}  to Lie
superalgebras and S. Benayadi, in \cite {Ben2000},  has used this
definition to study Lie superalgebras with reductive even part and
the action of the even part on the odd part is completely reducible
having these results been generalized in \cite {AB2004} for Malcev
superalgebras.

In case that we do not demand that the action of the even part on
the odd part is completely reducible we have a larger class of Lie
superalgebras  with reductive even part which contain the preceeding
ones and the remaining examples that appear in this study, obliged
us to use the notion of generalized double extension for Lie
superalgebras introduced in \cite{BBB} with the purpose to give an
inductive description of solvable quadratic  Lie superalgebras.

In this work we generalize all the results mentioned before,
presenting an inductive description of Malcev superalgebras with
reductive even part.  We prove that if $\displaystyle \frak {U}$ is
the set formed by $\{0\}$, the basic classical simple Lie
superalgebras, the non Lie simple Malcev algebra and the
one-dimensional Lie algebra, every Malcev superalgebra with
reductive even part is either an element of $\displaystyle \frak
{U}$ or is obtained from a finite number of elements of
$\displaystyle \frak {U}$ by a sequence of double extensions by the
one dimensional Lie algebra or/and generalized double extensions by
the one dimensional Lie superalgebra with null even part or/and
orthogonal direct sums.

\section{Basic notions}
For basic notions and properties concerning quadratic Malcev
superalgebras we propose to the reader a careful lecture of the
\cite{Alb1993,AB2004}. In this paper, we shall consider finite
dimensional Malcev superalgebras over an algebraically closed
commutative field $\displaystyle \mathbb{K}$ of characteristic zero.
\begin{defi} A superalgebra
$\displaystyle {M}={M}_{\bar 0} \oplus {M}_{\bar 1} $ (meaning a
$\displaystyle \mathbb{Z}_2$-graded algebra) is called a
\textit{Malcev superalgebra} if $\displaystyle XY=-(-1)^{xy} YX$,
$\displaystyle \forall_{ X \in {M}_x , Y \in {M}_y},$ and
\begin{eqnarray*}
\displaystyle
(-1)^{yz}(XZ)(YT)&=&((XY)Z)T+(-1)^{x(y+z+t)}((YZ)T)X\\
&&+(-1)^{(x+y)(z+t)}((ZT)X)Y+(-1)^{t(x+y+z )}((TX)Y)Z,
\end{eqnarray*}
for  $\displaystyle X\in {M}_{x}, Y \in {M}_{y}, Z \in {M}_{z}$, and
$\displaystyle T \in {M}_{t}$. We shall write $\displaystyle X \in
{M}_{x}$ to mean that $X$ is a homogeneous element of the Malcev
superalgebra $\displaystyle {M}$ of degree $x$, with $\displaystyle
x \in \mathbb{Z}_2$.
 \label{defi:defls}
\end{defi}
\begin{defi} Let $\displaystyle {M} $ be a  Malcev superalgebra and
$\displaystyle \phi : {M}\longrightarrow {M}$ an endomorphism of
$\displaystyle {M}$. We say that $\displaystyle \phi $ is a
\textit{Malcev operator} of $\displaystyle  {M}$ if
  \begin{eqnarray*}
\displaystyle
\phi((XY)Z)&=&(\phi(X)Y)Z-(-1)^{xy}\phi(Y)(XZ)-(-1)^{z(x+y)}(\phi(Z)X)Y\\
&&-(-1)^{x(y+z)}\phi(YZ)X,\hspace*{1cm}\forall_{X\in {M}_{x},Y\in {M}_{y},Z\in {M}_{z}}.
\end{eqnarray*}
We denote by $\displaystyle \Bigl(Op (M)\Bigr)_{\alpha}$ the vector
subspace of $\displaystyle End (M)$ formed by the Malcev operators
of $M$ of degree $\displaystyle \alpha$ ($\displaystyle \alpha \in
Z_2$). Then $\displaystyle Op (M)= \Bigl(Op (M)\Bigr)_{\bar
0}\oplus\Bigl(Op (M)\Bigr)_{\bar 1}$.
\end{defi}
\begin{defi}
Let $\displaystyle M$ be a Malcev superalgebra. A bilinear form  $B$
on $\displaystyle M$ is called
\begin{enumerate}
\item [(i)]
\textit{supersymmetric} \index{bilinear form (on a Lie
superalgebra)!supersymmetric} if $\displaystyle
  B(X,Y)=(-1)^{xy}B(Y,X)$, $\displaystyle \forall_{
X \in M_x , Y \in M_y}$ .
 \item [(ii)]  \textit{non-degenerate} \index{bilinear form (on a Lie
superalgebra)!non-degenerate} if $\displaystyle X \in M$ satisfies
$\displaystyle
  B(X,Y)=0$, $\displaystyle \forall_{
 Y \in M}$, then $\displaystyle X=0$.
 \item [(iii)]   \textit{invariant} \index{bilinear form (on a Lie
superalgebra)!invariant} if $\displaystyle
  B(XY,Z)=B(X,YZ)$, $\displaystyle \forall_{
X   , Y ,Z \in M}$.
 \item [(iv)]   \textit{even} \index{bilinear form (on a Lie
superalgebra)!even} if $\displaystyle
  B(X,Y)=0$, $\displaystyle \forall_{
 (X,Y)\in  M_\alpha \times M_\beta ,
  (\alpha ,\beta ) \in \mathbb{Z}_2 \times \mathbb{Z}_2 (   \alpha \neq \beta )}$.
\end{enumerate}
\end{defi}
\begin{defi}
A Malcev superalgebra $\displaystyle M$ is called \textit{quadratic}
if there
  exists a bilinear form $B$ on $\displaystyle M$ such that
  $B$ is even, supersymmetric, non-degenerate, and invariant. It is denoted by $\displaystyle (M,B)$ and
  $B$ is called an \textit{invariant scalar product} on $\displaystyle M$.
\end{defi}
\begin{defi}
Let $\displaystyle (M,B)$ be a quadratic Malcev superalgebra. A
homogeneous map $\phi \in End (M)$ of degree $\displaystyle \alpha$
(with $\displaystyle \alpha \in \mathbb{Z}_2$) is called
\textit{skew-supersymmetric}
 if
\begin{eqnarray*}
\displaystyle B(\phi(X),Y)&=&-(-1)^{\alpha x}B(X,\phi(Y)),
\hspace*{1cm}\forall_{  X \in M_x,  Y \in M}.
\end{eqnarray*}
We denote by $\displaystyle \Bigl(Op_a (M)\Bigr)_{\alpha}$ the
vector subspace of the skew-supersymmetric elements of
$\displaystyle \Bigl(Op(M)\Bigr)_{\alpha}$. We write  $\displaystyle
Op_a (M)= \Bigl(Op_a (M)\Bigr)_{\bar 0}\oplus\Bigl(Op_a
(M)\Bigr)_{\bar 1}$ which is a super-vector subspace of
$\displaystyle Op (M)$.
\end{defi}
\begin{defi} Let $\displaystyle {M}$  be a Malcev superalgebra and $V$ a $\displaystyle \mathbb{Z}_2$-graded vector space.
Let $\displaystyle \omega : {M}\times M\longrightarrow V$ be a
homogeneous bilinear map. If the assertion  $\displaystyle \omega
(X,Y)=-(-1)^{xy} \omega (Y,X)$, $\displaystyle \forall_{ X \in {M}_x
, Y \in {M}_y},$ and
  \begin{eqnarray*}
\displaystyle
(-1)^{yz}\omega(XZ,YT)&=&\omega((XY)Z,T)+(-1)^{x(y+z+t)}\omega((YZ)T,X)\\
&&+(-1)^{(x+y)(z+t)}\omega((ZT)X,Y)\\
&&+(-1)^{t(x+y+z )}\omega((TX)Y,Z),\hspace*{.5cm}
\forall_{ X\in {M}_{x}, Y \in {M}_{y}, Z \in {M}_{z},T \in {M}_{t}},
\end{eqnarray*}
 are satisfied we say that  $\displaystyle \omega$ is a \textit{Malcev $2$-cocycle} on
$\displaystyle {M}$ with values in $\displaystyle  V$.
\end{defi}
\begin{pro}[\cite{AB2004}]
Let $\displaystyle (M,B) $ be a quadratic Malcev superalgebra and
$\displaystyle \omega: M\times M \longrightarrow \mathbb{K}$ a
bilinear form of degree $\displaystyle \alpha \in \mathbb{Z}_2$.
\begin{enumerate}
  \item [(i)] There exists a homogeneous map $\phi \in  \mbox{End}(M)$ of
degree $\displaystyle \alpha$ such that
\begin{eqnarray*}
\displaystyle \omega (X,Y) &=&  B(\phi(X),Y), \hspace*{1cm}\forall_{X  , Y  \in M}. 
\end{eqnarray*}
  \item [(ii)] $\displaystyle \omega $
is a  Malcev $2$-cocycle on $\displaystyle M$ if and only if
$\displaystyle \phi$ is a skew-supersymmetric Malcev operator of
$M$.
\end{enumerate}
\end{pro}
\section{Generalized double extension of quadratic Malcev superalgebras}
We consider $\displaystyle \mathbb{K}e=(\mathbb{K}e)_{\bar 1} $ the
one-dimensional abelian Malcev superalgebra with even part zero and
$\displaystyle \mathbb{K}e^\ast$ its dual vector space.
\begin{pro}  Let $\displaystyle (M,B)$ be a quadratic Malcev
superalgebra and $\displaystyle D$ an odd skew-supersymmetric Malcev
operator of $\displaystyle (M,B)$. Let us consider the bilinear  map
$\displaystyle \varphi: M\times M \longrightarrow \mathbb{K}e^\ast$
defined by
\begin{eqnarray*}
\displaystyle \varphi (X,Y) &=&  -B(D(X),Y)e^\ast,
\hspace*{1cm}\forall_{
X , Y  \in M}. 
\end{eqnarray*}
Then $\displaystyle \varphi$ is a Malcev $2$-cocycle on $M$ with
values in $\displaystyle \mathbb{K}e^\ast$. Moreover, the
$\displaystyle \mathbb{Z}_2$-graded vector space $\displaystyle
M\oplus \mathbb{K}e^\ast$  with the multiplication defined by
\begin{eqnarray*}
\displaystyle (X+\alpha e^\ast)( Y+\beta
e^\ast)=XY+\varphi(X,Y),\hspace*{1cm}\forall_{(X+\alpha
e^\ast),(Y+\beta e^\ast) \in (M\oplus
\mathbb{K}e^\ast)}, 
\end{eqnarray*}
 is a Malcev superalgebra that is called the central
extension of $\displaystyle \mathbb{K}e^\ast$ by $\displaystyle M$
(by  means of $\displaystyle \varphi$).
\end{pro}
\begin{proof}
We will show that $\displaystyle M\oplus \mathbb{K}e^\ast$ is a
Malcev superalgebra. Since $D$ is odd and skew-supersymmetric we
have that $\displaystyle \varphi (X,Y) =  -(-1)^{xy}\varphi (X,Y)$,
for all $  X  \in M_x, Y \in M_y$. As the multiplication in $M$ is
graded skew-symmetric we conclude that also the multiplication
defined in $\displaystyle M\oplus \mathbb{K}e^\ast$ is graded
skew-symmetric. We have to ensure that
\begin{eqnarray*}
\displaystyle
&&(-1)^{yz}((X+\alpha e^\ast)(Z+\gamma e^\ast))((Y+\beta e^\ast)(T+\delta e^\ast))= \\
 &&\hspace*{1cm}=(((X+\alpha e^\ast)(Y+\beta e^\ast))(Z+\gamma e^\ast))(T+\delta e^\ast) \\
 &&\hspace*{1.5cm}+(-1)^{x(y+z+t)}(((Y+\beta e^\ast)(Z+\gamma e^\ast))(T+\delta e^\ast))(X+\alpha e^\ast) \\
 &&\hspace*{1.5cm}+(-1)^{(x+y)(z+t)}(((Z+\gamma e^\ast)(T+\delta e^\ast))(X+\alpha e^\ast))(Y+\beta e^\ast) \\
 &&\hspace*{1.5cm}+(-1)^{t(x+y+z )}(((T+\delta e^\ast)(X+\alpha e^\ast))(Y+\beta e^\ast))(Z+\gamma e^\ast) ,
\end{eqnarray*}
$\displaystyle \forall_{(X+\alpha e^\ast)\in (M\oplus
\mathbb{K}e^\ast)_{x},
 (Y+\beta e^\ast) \in (M\oplus \mathbb{K}e^\ast)_{y}, (Z+\gamma e^\ast) \in (M\oplus \mathbb{K}e^\ast)_{z} ,
 (T+\delta e^\ast) \in (M\oplus \mathbb{K}e^\ast)_{t}}$, or equivalently,
\begin{eqnarray*}
\displaystyle
&&(-1)^{yz}\bigl\{ (X Z ) (Y T ) + \varphi(X Z,Y T)  \bigr\}=\\
&& \hspace*{1cm}=( (X Y ) Z ) T +\varphi( (X Y ) Z,T)     \\
&& \hspace*{1.5cm}+(-1)^{x(y+z+t)}\bigl\{(( Y Z ) T ) X +\varphi(( Y Z ) T,X)  \bigr\} \\
 &&\hspace*{1.5cm}+(-1)^{(x+y)(z+t)}\bigl\{(( Z T ) X ) Y +\varphi(( Z T ) X,Y)   \bigr\} \\
&&\hspace*{1.5cm}+(-1)^{t(x+y+z )}\bigl\{(( T X ) Y ) Z +\varphi(( T X ) Y,Z)   \bigr\} .
\end{eqnarray*}
$\displaystyle \forall_{(X+\alpha e^\ast)\in (M\oplus
\mathbb{K}e^\ast)_{x},
 (Y+\beta e^\ast) \in (M\oplus \mathbb{K}e^\ast)_{y}, (Z+\gamma e^\ast) \in (M\oplus \mathbb{K}e^\ast)_{z} ,
 (T+\delta e^\ast) \in (M\oplus \mathbb{K}e^\ast)_{t}}$. Since $B$ is non-degenerate and $D$ is a Malcev operator then we have
\newline
$\displaystyle \forall_{(X+\alpha e^\ast)\in (M\oplus
\mathbb{K}e^\ast)_{x},
 (Y+\beta e^\ast) \in (M\oplus \mathbb{K}e^\ast)_{y}, (Z+\gamma e^\ast) \in (M\oplus \mathbb{K}e^\ast)_{z} ,
 (T+\delta e^\ast) \in (M\oplus \mathbb{K}e^\ast)_{t}}$,
\begin{eqnarray*}
\displaystyle
 (-1)^{yz} \varphi(X Z,Y T) &=& \varphi( (X Y ) Z,T) +(-1)^{x(y+z+t)} \varphi(( Y Z ) T,X)  \\
 && +(-1)^{(x+y)(z+t)} \varphi(( Z T ) X,Y) \\
 &&+(-1)^{t(x+y+z )} \varphi(( T X ) Y,Z).
\end{eqnarray*}
We infer that $\displaystyle \varphi$ is a Malcev $2$-cocycle of $M$
with values in $\displaystyle  \mathbb{K}e^\ast$. Using the fact
that $M$ is a Malcev superalgebra we conclude that $\displaystyle
M\oplus \mathbb{K}e^\ast$ is a Malcev superalgebra as required.
\end{proof}
\begin{pro}
Consider two Malcev superalgebras $\displaystyle M$ and
$\displaystyle V$, a linear map $\displaystyle \Omega:
M\longrightarrow Op (V)$, and an even skew-supersymmetric bilinear
map $\displaystyle \zeta: M\times M \longrightarrow V$ such that
\begin{enumerate}
  \item [(i)]$\displaystyle \forall_{X \in {M}_{x},
 Y  \in {M}_{y}, h \in {V}_{z} ,i \in V_{t}}$,
\begin{eqnarray*}
 &&\hspace*{0cm}  \bigl( \Omega(XY)(h)
\bigr)i- \Omega(X) \bigl(\Omega (Y)(h)i\bigr) -(-1)^{yz}\bigl(\Omega
(X)(h) \bigr)\bigl(\Omega (Y)(i) \bigr)\\
 && \hspace*{1.5cm}+(-1)^{xy} \Omega (Y)\bigl(\Omega (X)(hi) \bigr)
+(-1)^{tz+xy}\Omega (Y)\bigl(\Omega (X)(i)\bigr)
h\\
 && \hspace*{1.5cm}+\bigl(\zeta(X,Y)(h) \bigr)i   =  0;
\end{eqnarray*}
 \item [(ii)]$\displaystyle \forall_{X \in {M}_{x},
 Z  \in {M}_{z}, g \in {V}_{y} ,i \in V_{t}}$,
\begin{eqnarray*}
\displaystyle   &&\hspace*{-1cm} (-1)^{yz}\bigl\{ \zeta(X,Z)(gi)
+\Omega(XZ)(gi) \bigr\} =-(-1)^{z(x+y) }\Omega (Z)\bigl(\Omega
(X)(g)\bigr) i
\\
 && \hspace*{1.5cm} +(-1)^{yz}\Omega(X) \bigl(\Omega (Z)(g)i\bigr)
 -(-1)^{y(z+t) }\Omega (X)\bigl(\Omega (Z)(i)\bigr) g\\
 && \hspace*{1.5cm} +(-1)^{ty+(x+y)z}\Omega(Z) \bigl(\Omega
 (X)(i)g\bigr);
\end{eqnarray*}
  \item [(iii)] $\displaystyle \forall_{X \in {M}_{x},
 Y  \in {M}_{y}, T \in {M}_{t} ,h \in V_{z}}$,
\begin{eqnarray*}
  && \hspace*{0cm}
(-1)^{yz}\bigl(\Omega (X)(h) \bigr)\zeta(Y,T) +(-1)^{t(x+y+z) }
\Omega(T) \bigl(\zeta(X,Y)h\bigr) \\
 && \hspace*{1.5cm}+(-1)^{t(x+z) +xy} \Omega
(Y)\bigl(\zeta (T,X)\bigr) h  =  0;
\end{eqnarray*}
 \item [(iv)]$\displaystyle \forall_{X \in {M}_{x},
Y \in {M}_{y}, Z  \in {M}_{z} ,T \in M_{t}}$,
\begin{eqnarray*}
\displaystyle   &&\hspace*{-0.5cm} -(-1)^{yz}\Omega(XZ)
\bigl(\zeta(Y,T)\bigr)+(-1)^{t(y+z) }\Omega(X) \bigl(\Omega
(T)(\zeta(Y,Z))\bigr)\\
 && \hspace*{-0.5cm}+(-1)^{x(y+z+t) +yz}\Omega(Z) \bigl(\Omega
(Y)(\zeta(T,X))\bigr)- \Omega(X) \bigl( \zeta(YZ,T)
\bigr)\\
 && \hspace*{-0.5cm}-(-1)^{(x+y)(z+t)} \Omega(Z) \bigl( \zeta(TX,Y)
 \bigr)= \\
&&\hspace*{.2cm} =-(-1)^{yz}\bigl\{ \zeta(XZ,YT)
+\zeta(X,Z)\zeta(Y,T) \bigr\}
+\zeta((XY)Z,T)\\
 && \hspace*{.7cm}+(-1)^{x(y+z+t) }
\zeta((YZ)T,X)+(-1)^{(x+y)(z+t)}\zeta((ZT)X,Y)\\
 && \hspace*{.7cm}+(-1)^{t(x+y+z) }
\zeta((TX)Y,Z);
\end{eqnarray*}
  \item [(v)] $\displaystyle \forall_{X \in {M}_{x},
 Y  \in {M}_{y},Z \in M_{z}, i \in {V}_{t}} $,
\begin{eqnarray*}
\displaystyle  && \hspace*{0cm} (-1)^{yz}\Omega(XZ) \bigl(\Omega
(Y)(i)\bigr) -\zeta(XY,Z)i
-\Omega \bigl((XY)Z\bigr)(i) \\
 && \hspace*{1.5cm}+\Omega(X) \bigl(\Omega (YZ)(i)\bigr)
-(-1)^{xy} \Omega(Y) \bigl(\Omega
(X)(\Omega(Z)(i))\bigr)\\
 && \hspace*{1.5cm}+(-1)^{(x+y)z+xy} \Omega(Z) \bigl(\Omega
(Y)(\Omega(X)(i))\bigr)=  0.
\end{eqnarray*}
\end{enumerate}
 Then
the $\displaystyle \mathbb{Z}_2$-graded vector space $\displaystyle
M\oplus V$ endowed with the multiplication
\begin{eqnarray*}
\displaystyle (X+f)(Y+h)&=& (XY)_M + (fh)_V+\Omega(X)(h)-(-1)^{xy}
\Omega(Y)(f)+ \zeta(X,Y),
\end{eqnarray*}
where $\displaystyle (X+f) \in {(M\oplus V)}_{ x}$ and
$\displaystyle (Y+h) \in {(M\oplus V)}_{ y}$, is a Malcev
superalgebra that is called the \textit{generalized semi-direct
product} of $\displaystyle M$ and $\displaystyle V$  (by means of
$\displaystyle \Omega$ and $\displaystyle \zeta$). In particular, if
$\displaystyle \zeta=0$ then we have the \textit{semi-direct
product}  of $\displaystyle M$ and $\displaystyle V$  (by means of
$\displaystyle \Omega$) \cite{AB2004}. \label{pro:gsdpee}
\end{pro}
\begin{proof}
It is analogous to the proof of the \cite[Theorem 3.1]{AB2004}. To
show that $\displaystyle M\oplus V$ is a Malcev superalgebra it is
suffices to see the axiom of definition of a Malcev superalgebra
with $\displaystyle ((XY)Z)T$ in the following cases.
If all these elements are in $\displaystyle M$ it is clear since
$\displaystyle M $ is a Malcev superalgebra and we have assertion
(iv).
For the case where all these elements are in $V$ it is direct
because $V$ is a Malcev superalgebra.
In case where one of these elements is in $\displaystyle M$ and the
other three are in $\displaystyle V$ we know that $\displaystyle
\Omega (X)$ is a Malcev operator, for all $\displaystyle X \in M$.
The conditions (i) and (ii) imply that the axiom is valid in the
case $\displaystyle X \in {M}_{x},
 Y  \in {M}_{y}, h \in {V}_{z} ,i \in V_{t}$, and on the other hand, the assertions (iii) and (v)
  show the equality in the case
$\displaystyle X \in {M}_{x},
 Y  \in {M}_{y},Z \in M_{z}, i \in {V}_{t}$. Conversely, if $\displaystyle M\oplus V$ is a Malcev superalgebra
 these conditions are verified.
\end{proof}
Now, we shall use the particular case of generalized semi-direct
product of the Malcev superalgebras $\displaystyle M\oplus
\mathbb{K}e ^\ast$ and $\displaystyle \mathbb{K}e $ (where
$\displaystyle (M,B)$ is a quadratic Malcev superalgebra and
$\displaystyle (\mathbb{K}e)_{\bar 1} $ the one-dimensional abelian
Malcev superalgebra with even part zero)
\begin{teo}[]
Let $\displaystyle (M=M_{\bar 0}\oplus M_{\bar 1},B)$ be a quadratic
Malcev superalgebra. Suppose that $\displaystyle D \in (Op_a(M)
)_{\bar 1}$ and $\displaystyle A_{0} \in {M}_{{\bar 0}}$ such that
\begin{eqnarray}
\displaystyle && \hspace*{-1.7cm} D^2 (A_0) = \frac{1}{2}A_0A_0 \nonumber \\
&&  \hspace*{-1.7cm}  D(A_0 X )  = A_0D(X)-D(A_0)X,\label{eq:gdec1}\\
&&  \hspace*{-1.7cm} (A_0 X )Y  =  D(D(X)Y)+D^2(XY)+ (-1)^{x}D(X)D(Y)+(-1)^{xy}D^2(Y)X,\label{eq:gdec2}\\
&&  \hspace*{-1.7cm}  A_0 (XY)  =  D(D(X)Y)+D^2(X)Y -
(-1)^{xy}\bigl\{ D^2(Y)X+D( D(Y)X)\bigr\}\label{eq:gdec3},
\end{eqnarray}
where $\displaystyle X \in {M}_{x}, Y \in {M}_{y}$ are arbitrary.
Define a map $\displaystyle \Omega :(\mathbb{K}e)_{{\bar 1}}
\longrightarrow
 Op (M\oplus  \mathbb{K}e ^\ast) $ by $\displaystyle
\Omega (e)=\widetilde{D}$, where $\displaystyle
\widetilde{D}:M\oplus \mathbb{K}e ^\ast \longrightarrow M\oplus
\mathbb{K}e ^\ast $ satisfies $\displaystyle \widetilde{D}
(e^\ast)=0$ and
\begin{eqnarray*}
\displaystyle \widetilde{D}(X) &=& D(X) -(-1)^{ x} B(X,X_{0})e^\ast,
\hspace*{1cm}\forall_{
X \in M_x}. 
\end{eqnarray*}
Consider the map $\displaystyle \zeta :\mathbb{K}e \times
\mathbb{K}e \longrightarrow
 M\oplus  \mathbb{K}e ^\ast $ defined by $\displaystyle \zeta (e,e )=A_{0}$.
Then $\displaystyle {\frak k}=\mathbb{K}e \oplus M\oplus
 \mathbb{K}e  ^\ast$ equipped with the even
skew-symmetric bilinear map on $\displaystyle
 {\frak k}$ defined
by
\begin{eqnarray}
\displaystyle ee &=&  A_{0},\nonumber \\
\hspace*{0cm}eX  &=&  D(X)-(-1)^xB(X,A_{0}) e^\ast,
\hspace*{1cm}\forall_{
X \in {M}_x},\label{eq:bccnn}\\
\hspace*{0cm}XY &=& (XY)_{M}-B(D(X),Y)  e^\ast,
\hspace*{1cm}\forall_{
X,Y \in M},\nonumber \\
\hspace*{0cm} e^\ast {\frak k}  &=& \{0\},\nonumber
\end{eqnarray}
is the generalized semi-direct product of $\displaystyle M\oplus
\mathbb{K}e ^\ast$ by the one-dimensional Malcev superalgebra
$\displaystyle (\mathbb{K}e)_{{\bar 1}}$ (by means of $\displaystyle
\Omega$ and $\displaystyle \zeta$). Moreover the supersymmetric
bilinear form $\displaystyle \widetilde{B}: {\frak k}\times {\frak
k}\longrightarrow \mathbb{K}$ defined by
\begin{eqnarray*}
\displaystyle  && \widetilde{B}\mid_{M\times M}\: =\:B,\\
&& \widetilde{B}(e ,e^\ast)\: =\: 1,\\
&& \widetilde{B}(M,e)=  \widetilde{B}(M,e^\ast)\: =\: \{0\} ,\\
 &&
\widetilde{B}(e,e)=  \widetilde{B}(e^\ast,e^\ast)\: =\: 0 ,
\end{eqnarray*}
is an invariant scalar product on ${\frak k}$. The quadratic Malcev
superalgebra $\displaystyle ({\frak k},\widetilde{B}) $ is called
the generalized double extension  of $\displaystyle (M,B) $ by the
one-dimensional Malcev superalgebra $\displaystyle
(\mathbb{K}e)_{\bar 1}$ (by means of $D$ and $\displaystyle A_{0}
$).
 \label{teo:jjdufu}
\end{teo}
\begin{proof}
We start by showing that  $\displaystyle \widetilde{D}$ is a Malcev
operator of $\displaystyle M\oplus \mathbb{K}e ^\ast$. For all
$\displaystyle (X+\alpha e ^\ast) \in {(M\oplus \mathbb{K}e
^\ast)}_{ x}$, $\displaystyle (Y+\beta e ^\ast) \in {(M\oplus
\mathbb{K}e ^\ast)}_{ y}$, $\displaystyle (Z+\gamma e ^\ast) \in
{(M\oplus \mathbb{K}e ^\ast)}_{ z}$ we have to ensure that
\begin{eqnarray*}
\displaystyle &&\widetilde{D}\bigl(((X+\alpha e ^\ast)(Y+\beta e
^\ast))(Z+\gamma e ^\ast)\bigr)=
\bigl(\widetilde{D}(X+\alpha e ^\ast)(Y+\beta e ^\ast)\bigr)(Z+\gamma e ^\ast)\\
&&\hspace*{3cm}-(-1)^{xy}\widetilde{D}(Y+\beta e ^\ast)\bigl((X+\alpha e ^\ast)(Z+\gamma e ^\ast)\bigr)\\
&&\hspace*{3cm}-(-1)^{z(x+y)}\bigl(\widetilde{D}(Z+\gamma e ^\ast)(X+\alpha e ^\ast)\bigr)(Y+\beta e ^\ast)\\
&&\hspace*{3cm}-(-1)^{x(y+z)}\widetilde{D}\bigl((Y+\beta e ^\ast)(Z+\gamma e ^\ast)\bigr)(X+\alpha e ^\ast).
\end{eqnarray*}
Doing easy calculations, this condition is equivalent to
\begin{eqnarray*}
\displaystyle &&D ((X Y ) Z)
-(-1)^{ x+y+z}B((X Y ) Z, A_0)e ^\ast=\\
&&\hspace*{2cm}=(D(X) Y) Z
-B(D(D(X) Y),Z) e ^\ast \\
&&\hspace*{2.5cm}-(-1)^{xy}D(Y) (XZ)
+(-1)^{xy}B(D^2 (Y),XZ)e ^\ast\\
&&\hspace*{2.5cm}-(-1)^{z(x+y)}(D(Z) X) Y
+(-1)^{z(x+y)}B(D(D(Z) X),Y) e ^\ast\\
&&\hspace*{2.5cm}-(-1)^{x(y+z)}D (Y Z ) X
+(-1)^{x(y+z)}B(D^2(Y Z),X) e ^\ast.
\end{eqnarray*}
Since $D$ is a Malcev operator, it remains to see that
\begin{eqnarray*}
\displaystyle
 &&- B((X Y ) Z, A_0) =
-B(D(D(X) Y),Z)
+(-1)^{xy}B(D^2 (Y),XZ) \\
&&\hspace*{2cm}+(-1)^{z(x+y)}B(D(D(Z) X),Y)    +(-1)^{x(y+z)}B(D^2(Y Z),X) .
\end{eqnarray*}
As $D$ is skew-supersymmetric and $B$ is non-degenerate, it is   the
assertion \eqref{eq:gdec2}. Now we have to care about the several
conditions of Proposition \ref{pro:gsdpee}. Condition (v) is
immediate and (iv) comes direct from $\displaystyle   D^2 (A_0) =
\frac{1}{2}A_0A_0$. We have to verify that for all $\displaystyle
(Z+\gamma e ^\ast) \in {(M\oplus \mathbb{K}e ^\ast)}_{ z}$,
$\displaystyle (T+\delta e ^\ast) \in {(M\oplus \mathbb{K}e
^\ast)}_{ t}$,
\begin{eqnarray*}
 &&\hspace*{-1cm}  \bigl( \Omega(ee)(Z+\gamma e ^\ast)\bigr)(T+\delta e ^\ast)
 - \Omega(e) \bigl(\Omega (e)(Z+\gamma e ^\ast)(T+\delta e ^\ast)\bigr) \\
 && \hspace*{-1cm}-(-1)^{z}\bigl(\Omega (e)(Z+\gamma e ^\ast) \bigr)\bigl(\Omega (e)(T+\delta e ^\ast) \bigr)
- \Omega (e)\bigl(\Omega (e)((Z+\gamma e ^\ast)(T+\delta e ^\ast)) \bigr)\\
 && \hspace*{-1cm}-(-1)^{tz}\Omega (e)\bigl(\Omega (e)(T+\delta e ^\ast)\bigr)
(Z+\gamma e ^\ast) +\bigl(\zeta(e,e)(Z+\gamma e ^\ast)
\bigr)(T+\delta e ^\ast)  =  0,
\end{eqnarray*}
which we easily see that it is equivalent to
\begin{eqnarray*}
 &&\hspace*{-1cm}
-D(D(Z)T)- (-1)^{z}D(Z)D(T)-D^2(ZT)-(-1)^{zt}D^2(T)Z +(A_0 Z )T \\
&& -(-1)^{z+t}B(D(Z)T,A_0) e ^\ast -(-1)^{z+t}B(D(Z T),A_0) e ^\ast - B(D(A_0 Z ),T) e ^\ast\\
 && +(-1)^{z}B(D^2(Z),D(T)) e ^\ast+(-1)^{zt}B(D(D^2(T)), Z ) e ^\ast =0.
\end{eqnarray*}
 From \eqref{eq:gdec1} and \eqref{eq:gdec2}, since $B$ is even and
non-degenerate, and $D$ is an even skew-supersymmetric map, this
expression comes straightforward. Furthermore, we have to prove that
for all $\displaystyle (Y+\beta e ^\ast) \in {(M\oplus \mathbb{K}e
^\ast)}_{ y}$, $\displaystyle (T+\delta e ^\ast) \in {(M\oplus
\mathbb{K}e ^\ast)}_{ t}$,
\begin{eqnarray*}
\displaystyle   &&\hspace*{-1cm} (-1)^{y}\bigl\{ \zeta(e,e)((Y+\beta
e ^\ast)(T+\delta e ^\ast))
+\Omega(ee)((Y+\beta e ^\ast)(T+\delta e ^\ast)) \bigr\} \\
 && \hspace*{1.5cm}= (-1)^{y }\Omega (e)\bigl(\Omega
(e)(Y+\beta e ^\ast)\bigr) (T+\delta e ^\ast)\\
 && \hspace*{2cm} +(-1)^{y}\Omega(e) \bigl(\Omega (e)(Y+\beta e ^\ast)(T+\delta e ^\ast)\bigr)\\
  && \hspace*{2cm}-(-1)^{y(1+t) }\Omega (e)\bigl(\Omega (e)(T+\delta e ^\ast)\bigr) (Y+\beta e ^\ast)\\
 && \hspace*{2cm} -(-1)^{y(1+t) }\Omega(e) \bigl(\Omega
 (e)(T+\delta e ^\ast)(Y+\beta e ^\ast)\bigr),
\end{eqnarray*}
doing easy computation, it comes
\begin{eqnarray*}
\displaystyle   &&\hspace*{-1cm} A_0(Y T )-B(D(A_0), Y T) e ^\ast = D^2(Y) T-B( D^3(Y ),T) e ^\ast \\
 && \hspace*{2cm} +D(D(Y ) T )
 +(-1)^{y+t}B(D (Y) T,A_0) e ^\ast \\
  && \hspace*{2cm}-(-1)^{ yt} \bigl\{ D^2(T)Y -B(D^3(T),Y) e ^\ast\bigr\} \\
&& \hspace*{2cm}-(-1)^{ yt} \bigl\{ D(D(T)Y)+(-1)^{ y+t}B(D
(T)Y,A_0) e ^\ast \bigr\}.
\end{eqnarray*}
 From \eqref{eq:gdec1} and \eqref{eq:gdec3}, since $B$ is even and
non-degenerate, and $D$ is an even skew-supersymmetric map, we
obtain the assertion (ii) of the Proposition \ref{pro:gsdpee}.
Finally, we will show that for all $\displaystyle (Z+\gamma e ^\ast)
\in {(M\oplus \mathbb{K}e ^\ast)}_{ z}$,
\begin{eqnarray*}
  && \hspace*{0cm}
(-1)^{ z}\bigl(\Omega (e)(Z+\gamma e ^\ast) \bigr)\zeta(e,e) +(-1)^{
z }
\Omega(e) \bigl(\zeta(e,e)(Z+\gamma e ^\ast)\bigr) \\
 && \hspace*{1.5cm}+(-1)^{ z  } \Omega
(e)\bigl(\zeta (e,e)\bigr) (Z+\gamma e ^\ast) =  0,
\end{eqnarray*}
after some calculus, it is
\begin{eqnarray*}
  && \hspace*{0cm}
  D(Z)A_0+D(A_0Z)+D(A_0)Z\\
 && \hspace*{1.5cm} -B(D^2(Z),A_0) e ^\ast- (-1)^{ z  }B(A_0Z;A_0)e ^\ast-B(D^2(A_0),Z)e ^\ast  =  0.
\end{eqnarray*}
 From \eqref{eq:gdec1}, using $\displaystyle   D^2 (A_0) =
\frac{1}{2}A_0A_0$, since $B$ is even and non-degenerate, and $D$ is
skew-supersymmetric we get the desired expression. It is elementary
to see that
 $\displaystyle \mathbb{K}e \oplus M\oplus
 \mathbb{K}e ^\ast $ equipped with  multiplication defined by
\eqref{eq:bccnn} is the generalized semi-direct product of
$\displaystyle M\oplus \mathbb{K}e ^\ast $ by the one-dimensional
Malcev superalgebra $\displaystyle (\mathbb{K}e)_{{\bar 1}}$ (by
means of $\displaystyle \Omega$ and $\displaystyle \zeta$).
\end{proof}
We use the preceeding notion to present a class of quadratic Malcev
superalgebras with reductive even part and the action of the even
part on the odd part is not completely reducible
\begin{exa}
Let $\displaystyle ({M}={M}_{\bar 0}\oplus{M}_{\bar 1},B)$ be the
quadratic Malcev superalgebra such that dim${M}_{\bar 0}=1$ and the
dimension of ${M}_{\bar 1}$ is even. Let $\{a\}$ be a basis of
${M}_{\bar 0}$ and $\{v_1,v_2,...,v_n,y_1,y_2,...,y_n\}$ be a basis
of ${M}_{\bar 1}$, with multiplication defined by
\begin{eqnarray*}
 a.v_i=y_i=-v_i.a;\hspace*{0cm} M.y_i=\{0\},i=1,...,n; \\
v_i.v_j=\delta_{i,j}a, i,j=1,...,n;
\end{eqnarray*}
and invariant scalar product given by
\begin{eqnarray*}
\displaystyle B(a,a)=1;\hspace*{0cm}B(M_1,a)=\{0\},\\
B(y_i,y_j)=0=B(v_i,v_j);\hspace*{0cm}B(y_i,v_j)=\delta_{i,j}=-B(v_j,y_i), i,j=1,...,n;\\
\end{eqnarray*}

Consider $\{m_1,m_2,...m_n\}\subseteq K\setminus \{0\}$ and the
operator $D\in Op_a(M)_1$ defined by $D(a)=\sum_{i=1}^n m_iy_i;
D(v_i)=m_ia$ and $D(y_i)=0$. Then the generalized double extension
of (M,B) by the one dimensional Malcev superalagbera $Ke_1$ ( by
means of $D$ and a) is a quadratic Malcev superalgebra with
reductive even part and the action of the even part on the odd part
is not completely reducible.

In fact we have $\displaystyle {\frak k}=\mathbb{K}e \oplus M\oplus
 \mathbb{K}e  ^\ast$ equipped with the even
skew-symmetric bilinear map on $\displaystyle
 {\frak k}$ defined
by
\begin{eqnarray*}
ee = a;\hspace*{0cm}v_iv_j=\delta_{i,j}a; \hspace*{0cm}ev_i=m_ia; ae=\sum_{i=1}^n m_iy_i- e^\ast, \\
av_i= y_i-m_i e^\ast,\hspace*{0cm} e^\ast {\frak k}=y_i {\frak k}=
\{0\},
\end{eqnarray*}

Suppose now that $M_1$ is a completely reducible $M_0$-module. Then
there is an $M_0$-submodule $X$ such that $M_1=X\oplus Y$ where
$Y=<y_,...,y_n,e^\ast>$.But $M_0X\subseteq M_0M_1\subseteq Y$ and
$M_0X\subseteq X\cap Y=\{0\}$. As $M_0Y=\{0\}$  then $M_0M_1=\{0\}$
which is a contradiction.
\end{exa}

Now we shall show the converse of the Theorem  \ref{teo:jjdufu}
\begin{teo}
Let $\displaystyle ({M}={M}_{\bar 0}\oplus{M}_{\bar 1},B)$ be a
$B$-irreducible quadratic Malcev superalgebra such that
$\displaystyle \dim M > 1$. If $\displaystyle \frak z(M)\cap
{M}_{\bar 1}\neq \{0\}$ then $\displaystyle ({M},B) $ is a
generalized double extension of a quadratic Malcev superalgebra
$\displaystyle (N,\widetilde{B}) $ such that $\displaystyle \dim N =
\dim M -2$  by the one-dimensional Malcev
 superalgebra with even part zero.
 \label{teo:2345}
\end{teo}
We will prove the result by showing the following sequence of claims
as it is done in  case of the quadratic Lie superalgebras
\cite{ABB1}, and in case of the odd-quadratic Lie superalgebras
\cite{ABB2}. First, we will determine the quadratic Malcev
superalgebra $\displaystyle (N,\widetilde{B}) $; then we will show
that the
 quadratic Malcev superalgebra $\displaystyle (M,B) $
is the generalized double extension of $\displaystyle
(N,\widetilde{B}) $ by the one-dimensional Malcev superalgebra
$\displaystyle (\mathbb{K}e)_{{\bar 1}}$.
\begin{proof} Let us
assume that $\displaystyle (M,B) $ is a $B$-irreducible quadratic
Malcev superalgebra such that $\displaystyle \dim M > 1$ and
$\displaystyle \frak z(M)\cap M_{\bar 1}\neq \{0\}$. We set
$\displaystyle e^\ast $ a non-zero element of $\displaystyle \frak
z(M)\cap M_{\bar 1}$ and denote $\displaystyle I=\mathbb{K}e^\ast$.
As $B$ is even we have $\displaystyle M_{\bar 0}\subseteq J$, where
$J$ is the orthogonal of $I$ with respect to $B$. Since $B$ is
non-degenerate and even then there exists $\displaystyle e \in
M_{\bar 1} $ such that $\displaystyle B(e^\ast,e)\neq 0$. We may
assume that $\displaystyle B(e^\ast,e)=1$. As $\displaystyle e
\notin J$ and $\displaystyle \dim J = \dim M -1$ we infer that
$\displaystyle  M =J \oplus \mathbb{K}e$. Consider  the
two-dimensional vector subspace $\displaystyle  A = \mathbb{K}e^\ast
\oplus \mathbb{K}e$ of $\displaystyle M$. Since $\displaystyle
B\mid_{A\times A}$ is non-degenerate we have $\displaystyle M= A
\oplus N $, where $\displaystyle N $ is the orthogonal of $A$ with
respect to $B$. It comes that $\displaystyle \widetilde{B}
=B\mid_{N\times N}$ is non-degenerate. As $B$ is even we have
$\displaystyle \mathbb{K}e^\ast \subseteq J$, and so $\displaystyle
\mathbb{K}e^\ast\oplus N \subseteq J$. From $\displaystyle \dim \:(
\mathbb{K}e^\ast\oplus N) = \dim M -1= \dim J$ it comes that
$\displaystyle J=\mathbb{K}e^\ast \oplus N$. So $\displaystyle N$ is
a graded vector subspace of $\displaystyle M$ contained in the
graded ideal $\displaystyle J=N\oplus \mathbb{K}e^\ast$ of
$\displaystyle M$. Then we have
\begin{eqnarray*}
\displaystyle XY &=& \alpha (X,Y) + \varphi (X,Y)e^\ast,
\hspace*{1cm}\forall_{ X , Y \in N}, 
\end{eqnarray*}
where $\displaystyle \alpha (X,Y) \in N$ and $\displaystyle
 \varphi (X,Y) \in \mathbb{K}$. Further,
\begin{eqnarray*}
\displaystyle eX &=& D (X) + \psi(X)e^\ast,
\hspace*{1cm}\forall_{ X \in N},
\end{eqnarray*}
where $\displaystyle D (X) \in N$ and $\displaystyle
 \psi (X) \in \mathbb{K}$.\\
 \\
\textbf{Claim 1.} Then  $\displaystyle
 N  $ endowed with the multiplication $\displaystyle \alpha$ is a Malcev superalgebra and
 $\displaystyle \widetilde{B} $ is an invariant scalar product on $\displaystyle
 N$.\\
 \\
\textit{ Proof of Claim 1}: From graded skew-symmetry on
$\displaystyle M$ we get that
\begin{eqnarray}
\displaystyle \alpha (X,Y)&=&-(-1)^{xy} \alpha (Y,X),
\hspace*{1cm}\forall_{ X \in N_x , Y \in N_y},
\label{eq:ff1}\\
 \varphi (X,Y)&=&-(-1)^{xy} \varphi (Y,X),
\hspace*{1cm}\forall_{ X \in N_x , Y \in N_y}.\nonumber
\end{eqnarray}
Moreover, observing that for all $\displaystyle X \in N_x,Y \in
N_y,Z \in N_z,T \in N_t$
\begin{eqnarray*}
\displaystyle
 (XZ)(YT)&=&\alpha(\alpha(X,Z),\alpha(Y,T))+\varphi( \alpha(X,z),\alpha(Y,T))e^\ast,
\end{eqnarray*}
\begin{eqnarray*}
\displaystyle
((XY)Z)T&=&\alpha(\alpha(\alpha(X,Y)Z)T)+\varphi(\alpha(\alpha(X,Y),Z),T)e^\ast,
\end{eqnarray*}
by the  second property of definition of  Malcev superalgebras, it
comes that for all $\displaystyle X \in N_x,Y \in N_y,Z \in N_z,T
\in N_t$
\begin{eqnarray}
\displaystyle
(-1)^{yz}\alpha(\alpha(X,Z),\alpha(Y,T))&=&\alpha(\alpha(\alpha(X,Y),Z),T) \nonumber  \\
&&+(-1)^{x(y+z+t)}\alpha(\alpha(\alpha(Y,Z),T),X) \nonumber  \\
&&+(-1)^{(x+y)(z+t)}\alpha(\alpha(\alpha(Z,T),X),Y) \nonumber  \\
&&+(-1)^{t(x+y+z
)}\alpha(\alpha(\alpha(T,X),Y),Z),\hspace*{1cm}\label{eq:ff3}
\end{eqnarray}
and
\begin{eqnarray*}
\displaystyle
(-1)^{yz}\varphi(\alpha(X,Z),\alpha(Y,T))&=&\varphi(\alpha(\alpha(X,Y),Z),T)\\
&&+(-1)^{x(y+z+t)}\varphi(\alpha(\alpha(Y,Z),T),X)\\
&&+(-1)^{(x+y)(z+t)}\varphi(\alpha(\alpha(Z,T),X),Y)\\
&&+(-1)^{t(x+y+z )}\varphi(\alpha(\alpha(T,X),Y),Z).
\end{eqnarray*}
 By \eqref{eq:ff1} and
\eqref{eq:ff3} we infer that $\displaystyle N $ endowed with the
multiplication $\displaystyle \alpha $ is a Malcev superalgebra and
we denote $\displaystyle \alpha= (,)_N $. It is immediate that
$\widetilde{B}$ is an
 invariant scalar product on $\displaystyle N$ completing the proof of Claim
 1.\\
\\
Since $\displaystyle e \in M_{\bar 1} $ then $\displaystyle [e,e] $
is not necessarily zero, and we can write
\begin{eqnarray*}
\displaystyle [e,e]&=&  A_0 ,
\end{eqnarray*}
where $\displaystyle A_0 \in N_{\bar 0} $.\\
\\
\textbf{Claim 2.}  Then $\displaystyle D$ is an odd
skew-supersymmetric Malcev operator of $\displaystyle
(N,\widetilde{B})$  such that
\begin{eqnarray}
\displaystyle && \hspace*{-1.7cm}  D^2 (A_0) = \frac{1}{2}(A_0A_0)_N \nonumber \\
&&  \hspace*{-1.7cm}D(A_0 X )  = A_0D(X)-D(A_0)X,\label{eq:gdec1a}\\
&&  \hspace*{-1.7cm} (A_0 X )Y  =  D(D(X)Y)+D^2(XY)+ (-1)^{x}D(X)D(Y)+(-1)^{xy}D^2(Y)X,\label{eq:gdec2a}\\
&&  \hspace*{-1.7cm}  A_0 (XY)  =  D(D(X)Y)+D^2(X)Y -
(-1)^{xy}\bigl\{ D^2(Y)X+D( D(Y)X)\bigr\}\label{eq:gdec3a},
\end{eqnarray}
for every $\displaystyle X \in {M}_{x}, Y \in {M}_{y}$. Moreover,
$\displaystyle (M,B)$ is the generalized double extension   of
$\displaystyle (N,\widetilde{B}) $ by the one-dimensional Malcev
superalgebra $\displaystyle (\mathbb{K}e)_{{\bar 1}}$ (by means of
$D$ and $\displaystyle A_{0} $).\\
\\
\textit{ Proof of Claim 2}: We start by proving that $\displaystyle
D$ is an odd skew-supersymmetric Malcev operator of $\displaystyle
(N,\widetilde{B})$. It is clear that $\displaystyle D$ is a
homogeneous map of degree $\displaystyle {\bar 1}$. The second
property of definition of Malcev superalgebras for all
$\displaystyle X \in N_x,Y \in N_y,Z \in N_z$
\begin{eqnarray*}
\displaystyle
(-1)^{yz}(XZ)(Ye)&=&((XY)Z)e+(-1)^{x(y+z+1)}((YZ)e)X\\
&&+(-1)^{(x+y)(z+1)}((Ze)X)Y+(-1)^{  x+y+z }((eX)Y)Z,
\end{eqnarray*}
implies that  $\displaystyle D$ is a
 Malcev operator of $\displaystyle N$. On the other hand, Using  the invariance of
$B$, $\displaystyle B([e,X],Y)=-(-1)^{x}B(X,[e,Y])$, for every
$\displaystyle  X \in N_x,  Y \in N $, we obtain
\begin{eqnarray*}
\displaystyle \widetilde{B}\bigl(D(X),Y\bigr)&=&-(-1)^{
x}\widetilde{B}\bigl(X,D(Y)\bigr), \hspace*{1cm}\forall_{  X \in
N_x,  Y \in N},
\end{eqnarray*}
which means that $\displaystyle D \in  \Bigl(Op_a
(N,\widetilde{B})\Bigr)_{\bar 1}$. Due to $\displaystyle
-(ee)(ee)=-2((ee)e)e $ we obtain $\displaystyle D^2 (A_0) =
\frac{1}{2}(A_0A_0)_N$.
 From $\displaystyle B([e,X],Y)=B(e,[X,Y])$, for all $\displaystyle   X ,
Y \in N $, we get that
\begin{eqnarray*}
\displaystyle \varphi(X,Y) &=& -B(D(X),Y), \hspace*{1cm}\forall_{ X
, Y \in N}.
\end{eqnarray*}
Due to $\displaystyle B([e,X],e)= B(e,[X,e])$, for all
$\displaystyle X \in N_x $, we conclude that
\begin{eqnarray*}
\displaystyle \psi(X)&=&-(-1)^{ x}B(X,A_0), \hspace*{1cm}\forall_{ X
\in N_x}.
\end{eqnarray*}
Using sucessively the second property of definition of Malcev
superalgebras
\begin{eqnarray*}
\displaystyle
-(ee)(eX)&=&((ee)e)X+(-1)^{x }((ee)X)e + ((eX)e)e\\
&&+(-1)^{ x }((Xe)e)e,\\
(-1)^{x}(eX)(eY)&=&((ee)X)Y-(-1)^{ x+y }((eX)Y)e+((XY)e)e\\
&&+(-1)^{yx}((Ye)e)X,\\
(-1)^{x }(ee)(XY)&=&((eX)e)Y-(-1)^{ x+ y)}((Xe)Y)e \\
&&+(-1)^{(1+x)(1+y)}((eY)e)X+(-1)^{yx}((Ye)X)e,
\end{eqnarray*}
for  $\displaystyle X\in {M}_{x}, Y \in {M}_{y}$, we infer
respectively assertions \eqref{eq:gdec1a}, \eqref{eq:gdec2a}, and
\eqref{eq:gdec3a}, which states Claim 2.\\
\\
Using the results of Claims 1 and 2 we conclude Theorem
\ref{teo:2345}.
\end{proof}

\section{Structure and inductive classification of quadratic Malcev superalgebras with reductive even part}
\begin{teo}(proved in \cite{AB2004})
Let $\displaystyle ({M}={M}_{\bar 0}\oplus{M}_{\bar 1},B)$ be a
$B$-irreducible quadratic Malcev superalgebra such that
$\displaystyle \dim M > 1$. If $\displaystyle \frak z(M)\cap
{M}_{\bar 0}\neq \{0\}$ then $\displaystyle ({M},B) $ is a double
extension of a quadratic Malcev superalgebra $\displaystyle
(N,\widetilde{B}) $ such that $\displaystyle \dim N = \dim M -2$  by
the one-dimensional Lie
 algebra.
 \label{teo:234}
\end{teo}
So to conclude our study we must analyse the case where
$\displaystyle M$ is a $B$-irreducible quadratic Malcev superalgebra
with reductive even part  such that $\displaystyle \frak
z(M)=\{0\}$. Analogously to the Lie case (cf. \cite{ABB1}), the
socle of $\displaystyle M$ has a remarkable importance in this study
and we can prove with the same reasoning that if $\displaystyle M$
is a $B$-irreducible quadratic Malcev superalgebra with reductive
even part that is neither the one dimensional Lie algebra nor a
simple Malcev superalgebra the socle coincides with the center. As a
consequence of this result, with a proof analogously to the one in
Lie case \cite{ABB1}, we characterize quadratic Malcev superalgebras
with reductive even part that have null center as semisimple Malcev
superalgebras.

Finally we present an inductive description of quadratic Malcev
superalgebras such that the even part is a reductive Malcev algebra.
Clearly we improve the results proved in \cite{AB2004} since we work
with a larger class of superalgebras.

Let $\displaystyle \frak {U}$ be the set formed by $\displaystyle \{
0\}$, basic classical Lie superalgebras, simple (non Lie) Malcev
algebra and one-dimensional Lie algebra.
\begin{teo}[]
Let $\displaystyle ({M}={M}_{\bar 0}\oplus{M}_{\bar 1},B) $ be a
quadratic Malcev superalgebra such that $\displaystyle {M}_{\bar 0}
$ is a reductive Malcev algebra. Then either $\displaystyle M$ is an
element of $\displaystyle \frak U$, or $\displaystyle M$ is obtained
by a sequence of double extensions by the one-dimensional Lie
algebra, or generalized double extensions by the one-dimensional
Malcev superalgebra with even part zero, and/or by orthogonal direct
sums of quadratic Malcev superalgebras from a finite number of
elements of $\displaystyle \frak U$.\label{teo:ccc}
\end{teo}
\begin{proof} The proof is by induction on the dimension of $\displaystyle M$ and is analogous to the corresponding theorem
in Lie case (cf. \cite{ABB1}), using the fact that  a simple
quadratic Malcev superalgebra is either a basic classical Lie
superalgebra or is the simple (non Lie) Malcev algebra.
\end{proof}

\end{document}